\documentclass[10pt,leqno]{amsart}
\usepackage{txfonts}
\usepackage{amsmath,amsthm,amsfonts}
\usepackage{latexsym}
\usepackage{amssymb}
\usepackage{mathrsfs}
\usepackage[colorlinks,
            linkcolor=red,
            anchorcolor=blue,
            citecolor=green,
            ]{hyperref}
\usepackage[dvips]{epsfig}

\newcounter{RomanNumber}

\renewcommand{\thesection}{\arabic{section}}

\newtheorem{theorem}{Theorem}[]
\newtheorem{lemma}{Lemma}[section]

\newtheorem{defi}{Definition}[section]

\theoremstyle{remark}
\newtheorem{remark}{Remark}[section]

\renewcommand{\theequation}{\thesection .\arabic{equation}}
\let\sect\section
\renewcommand\section{\setcounter{equation}{0}
\gdef\theequation{\thesection .\arabic{equation}}\sect}

\makeatletter

\newcommand{\Rmnum}[1]{\expandafter\@slowromancap\romannumeral #1@}
\makeatother

\newcommand{\cx}{{\mathcal{X}}}

\newcommand{\IC}{{\mathbb{C}}}

\newcommand{\IR}{{\mathbb{R}}}

\newcommand{\be}{\begin{eqnarray}}
\newcommand{\ee}{\end{eqnarray}}

\newcommand{\me}{\mathop{\rm{meas}\, }}

\newcommand{\un}{\underline}

\def\beeq{\begin{equation}}
\def\eneq{\end{equation}}

\def\bm{\begin{matrix}}
\def\endm{\end{matrix}}

\def\Im{{\rm Im}}

\begin{document}

\title[Non-perturbative positive Lyapunov exponent and  applications]{Non-perturbative positive Lyapunov exponent of Schr\"odinger  equations and  its applications to the skew-shift }

\author{Kai Tao}

\address{College of Sciences, Hohai University, 1 Xikang Road Nanjing Jiangsu 210098 P.R.China }

\email{ktao@hhu.edu.cn,\ tao.nju@gmail.com}

\thanks{The author was supported by  the National Nature Science Foundation of China (Grant 11401166) and the Fundamental Research Funds for the Central Universities(Grant 2017B17214). He also wishes to thank Jiangong You for helpful discussions. }

\date{}

\begin{abstract}We first study  the  discrete Schr\"odinger  equations with  analytic potentials given by a class of transformations. It is shown that if the coupling number is large, then its logarithm equals approximately to the Lyapunov exponents. When the transformation becomes the  skew-shift, we prove  that the Lyapunov exponent  is  week H\"older continuous, and the spectrum satisfies Anderson Localization and contains large intervals. Moreover, all of these  conclusions  are non-perturbative.
\end{abstract}

 \maketitle

\section{Introduction}
Denote by $(\mathbb{Y},\mathscr{B},m)$ the probability space and $\mathbb{T}:=\mathbb{R}/\mathbb{Z}$ the torus equipped with its Haar measure. Let the measure preserving transformation $T:\mathbb{T}\times \mathbb{Y}\to \mathbb{T}\times \mathbb{Y}$ has the form
\begin{equation}
\label{map} T(x,y)=\left (x+f(y),g(y)\right ),
\end{equation}
with $g:\mathbb{Y}\to \mathbb{Y}$ and $f:\mathbb{Y}\to\mathbb{T}$.

In this paper, we first consider the following  discrete Schr\"odinger equations on $l^2(\mathbb{Z}^+)$:
\begin{equation}\label{schequ1}
( S_{(x,y),\lambda v}\phi)(n)=\phi(n+1)+\phi(n-1)+\lambda v(\pi_{\mathbb{T}}(T^n (x,y)))\phi(n)=E\phi(n),
\end{equation}
where  $v(x)$ is a real analytic function on $\mathbb{T}$, $\pi_{\mathbb{T}}$ is a projection from $\mathbb{Y}\times \mathbb{T}$ to $\mathbb{T}$. Then the equations  (\ref{schequ1}) can be expressed as
\[
\left (\begin{array}{cc}
  \phi(n+1) \\ \phi(n) \\
\end{array} \right )=M_n(x,y,E)\left (\begin{array}{cc}
  \phi(1) \\ \phi(0) \\
\end{array} \right ),
\]where
\[ M_n(x,y,E)=\prod_{j=n}^1\left ( \begin{array}{cc}
 E-\lambda v(\pi_{\mathbb{T}}(T^j(x,y))) & -1 \\
 1& 0 \\
  \end{array}\right )\]is called the transfer matrix of (\ref{schequ1}). Define
  \[ L_n(E)=\frac{1}{n}\int_{\mathbb{Y}}\int_{\mathbb{T}}\log \| M_n(x,y,  E)\|dxdm(y),\]and by the subadditive property, the limit
  \begin{equation}\label{lyaexp}
L(E)= \lim\limits_{n\to \infty} L_n(E)
\end{equation}exists, which is called the Lyapunov exponent of (\ref{schequ1}). Note that $\det M_n=1$, which implies $L(E)\ge 0$. But, in large coupling regimes, it is always positive as follow:
\begin{theorem}\label{positive}
  For any $\kappa>0$, there exists $\lambda_0=\lambda_0(v,\kappa)$ such that if $E$ is in the spectrum of (\ref{schequ1}) and $\lambda>\lambda_0$, then
  \[(1-\kappa)\log\lambda<L(E)<(1+\kappa)\log\lambda.\]
\end{theorem}
\begin{remark}\label{re}
  Due to the uniform hyperbolicity, the Lyapunov exponent is always positive when $E$ is regular.
\end{remark}

The transformation (\ref{map}) has many instances. The most famous one, called the shew-shift mapping, is defined on $\mathbb{T}^d$:
\begin{equation}
\label{skewshift}T_{d,\omega}(\un x=(x_1,\cdots,x_d))=(x_1+x_2,x_2+x_3,\cdots,x_{d-1}+x_d,x_d+\omega).
\end{equation}Then,  the equations (\ref{schequ1}) become
\begin{equation}\label{schequ}
 ( S_{\un x,\omega,\lambda v}\phi)(n)=\phi(n+1)+\phi(n-1)+\lambda v(\pi_1(T^n_{d,\omega} (\un x)))\phi(n)=E\phi(n),
\end{equation}
where $\pi_j$ is a projection from $\mathbb{T}^d$ to its $j$-th coordinate. This transformation $T_{d,\omega}$ is ergodic on $\mathbb{T}^d$ with any irrational $\omega$. Thus, due to the Kingman's Subadditive Ergodic Theorem, we have
\[ L(E,\omega)=\lim\limits_{n\to \infty}\frac{1}{n}\log \|M_n(\un x,E,\omega)\|\ \mbox{for almost every } \un x\in\mathbb{T}^d.\]

In this paper, when we say that $\omega\in(0,1)$ is the
Diophantine number(DN), it means that $\omega$ satisfies the Diophantine condition
\begin{equation}\label{13001}
\|n\omega\| \geq \frac{C_{\omega}}{n(\log n)^\alpha}\ \ \mbox{for
all}\ n\not=0.\end{equation} It is well known that for a fixed $\alpha>1$
almost every $\omega$ satisfies (\ref{13001}).

Now we have the following theorem for the Schr\"odinger equations (\ref{schequ}):
\begin{theorem}\label{skt}
 Assume that $v(x)$ is a real analytic function on $\mathbb{T}$. Then there exists $\lambda_0=\lambda_0(v)>0$ such that the followings hold for any $d\ge 1$:
\begin{enumerate}
  \item[{\rm{(P)}}]For any irrational $\omega$, if the coupling number $\lambda>\lambda_0$, then the Lyapunov exponent $L(E,\omega)$ of (\ref{schequ}) is positive for all $E\in\mathbb{R}$:
  \[L(E,\omega)\ge \frac{99}{100}\log \lambda >0.\]
  \item[{\rm{(C)}}]Let $\omega$ be the Diophantine number and $\lambda>\lambda_0$. Then  $L(E,\omega)$ is a continuous function of E with modulus of continuity
      \[
      h(t)=\exp\left (-c|\log t|^\tau\right ),
      \]
      where  $\tau$ is an absolute positive constant and $c=c(v)$ is also positive.
  \item[{\rm{(AL)}}]Let $\lambda>\lambda_0$ and fix $\un x\in\mathbb{T}^d$. Then for almost every $\omega\in\mathbb{R}$, the Schr\"odinger operator $ S_{\un x,\omega,\lambda v}$ satisfies Anderson localization, i.e., it  has pure point spectrum with exponentially decaying eigenfunctions.
  \item[{\rm{(I)}}]Let $d> 1$, $\omega$ be the Diophantine number and $\sigma( S_{\un x,\omega,\lambda v})$ denote the spectrum of the operator $ S_{\un x,\omega,\lambda v}$. For any given $\delta>0$, there exists $\lambda'_0=\lambda'_0(\delta, v)>0$ such that for any $\lambda>\max\{\lambda_0,\lambda'_0\}$,
      \[\lambda\mathcal{E}_\delta\subset \sigma( S_{\un x,\omega,\lambda v}),\ \forall x\in \mathbb{T}^d,\]
  where
  \[\mathcal{E}_\delta=\{E|\  \exists  x: v(x)=E\ \mbox{and }|v'(x)|\ge\delta\}\]is a union of intervals.
\end{enumerate}
\end{theorem}
\begin{remark}
  \begin{enumerate}
    \item[{\rm{(1)}}]In \cite{BGS}, Bourgain, Goldstein and Schlag studied the following shew-shift Schr\"odinger equations:
    \begin{equation}\label{gensk}
      (H_{\un x,\omega,\lambda v}\phi)(n)=\phi(n+1)+\phi(n-1)+\lambda v(T^n_{d,\omega} (\un x))\phi(n)=E\phi(n),
    \end{equation}
    where $v(x)$ is a real analytic function on $\mathbb{T}^d$. They proved the positive Lyapunov exponent, week H\"older continuity and Anderson localization, too. But all of these conclusions are perturbative and the set of the suitable $\omega$ only has positive measure. ~\\

    \item[{\rm{(2)}}]Obviously, our model is a special case of the one in \cite{BGS}. But in fact, people always pay more attention to ours. For example, let $d=2$ and $v(x)=\cos x$, then the Schr\"odinger equations (\ref{schequ})  become  the following Almost Mathieu equations
       \[\left(M_{(x,y),\omega,\lambda v}\phi\right)(n)=\phi(n+1)+\phi(n-1)+\lambda \cos\left  (x+ny+\frac{n(n-1)}{2}\omega\right )\phi(n)=E\phi(n).\]
    Bourgain  conjectured in \cite{B} that
    \begin{enumerate}
    \item[{\rm{(i)}}]If $\lambda\not=0$, then $L(E)$ is positive for all energies;
    \item[{\rm{(ii)}}] For all $\lambda\not=0$ and $(x,y)\in\mathbb{T}^2$, the  operator $M_{(x,y),\omega,\lambda v}$  has pure point spectrum with Anderson localization;
    \item[{\rm{(iii)}}]There are no gaps in the spectrum.
    \end{enumerate}Thus, \cite{BGS} gave the perturbative results of (i) and (ii) in the large coupling regimes.~\\

    \item[{\rm{(3)}}]The first result of (iii) is \cite{K2}. In this reference, Kr\"uger  developed the theory of parameterizing isolated eigenvalues and applied the perturbative Large Deviation Theorem(LDT for short) from \cite{BGS} to show that the spectrum of (\ref{schequ}) has intervals with large $\lambda$. Thus, we can improve Kr\"uger's conclusion to our Statement (I),  because we prove the non-perturbative LDT in Section 3.~\\

    \item[{\rm{(4)}}]So, compared to \cite{BGS} and \cite{K2}, one of the highlights of our paper is that it is the first to give  the non-perturbative answers to Bourgain's conjecture for the skew-shift Schr\"odinger equations in the large coupling regimes. Meanwhile, this non-perturbation makes the set of our suitable $\omega$ have full measure. Note that, when $d=1$, (\ref{schequ}) becomes the quasi-periodic Schr\"odinger equations  and their non-perturbative Anderson localization had been proven in \cite{BG}. Thus, the third highlight is that our Statement (AL) generalizes it for  $d>1$.~\\

    \item[{\rm{(5)}}]There are also some other  works by Kr\"uger concerning about the positive of the  Lyapunov exponent of our model. In \cite{K1},  he proved that  there exists $c>0$ such that
    \[ meas \{E:L(E)<c\}\to 0\] when $d\to +\infty$; what's more, for any $\epsilon>0$, there exists $\lambda_1(d,\epsilon)$ such that
    \[meas \{E:L(E)<\log \lambda\}<\epsilon\]for all $\lambda>\lambda_1$. Further, if  $\omega$ is a Diophantine number and $v$ is  a trigonometric polynomial of degree $K$, then there exists $\lambda_2(K,d,\omega)$ such that
    \[
      L(E)>\frac{1}{100}\log \lambda
    \]
for all $E$, when $\lambda>\lambda_2$. In \cite{K4}, for the non-degenerate potential $v(x)$, i.e., there exist $F$ and $\alpha$ such that for any $E\in\mathbb{R}$ and $\epsilon>0$,
\[meas\{x:|v(x)-E|<\epsilon\}\leq F\epsilon^\alpha,\]he proved that there exist $\lambda_3=\lambda_3(v)$ and $\kappa=\kappa(v)$ such that for $\lambda>\lambda_3$, there exists a set $\mathcal{E}_{\alpha ,\lambda}$ of measure
\[|\mathcal{E}_{\alpha,\lambda}|\leq \exp(-\lambda^{\frac{\alpha}{2}})\]
such that for $E\not\in \mathcal{E}_{\alpha,\lambda}$,
\[L(E)>\kappa \log \lambda.\]  Above all, it is easily seen  that our Statement (P)  is  optimal.~\\

\item[{\rm{(6)}}]When $d=1$, Goldstein and Schlag shew that the Lyapunov exponent is H\"older continuous of E in \cite{GS}. It may be right for $d>1$. But until now, we have no idea to get it.
  \end{enumerate}
\end{remark}~\\

It is obvious that  Statement (P) in Theorem \ref{skt} comes directly from Theorem \ref{positive} with $\kappa=\frac{1}{100}$. So, we organize this paper as follows. In Section 2, we develop Bourgain and Goldstein's method, which was applied to the quasi-periodic  Schr\"odinger equations  in \cite{GS}, to prove Theorem \ref{positive}. Then  the most important lemma, the non-perturbative Large Deviation Theorem for (\ref{schequ}), is given in Section 3. Finally, the proofs of Statements (C), (AP) and (I) are presented in the last section.

\section{positive Lyapunov exponent}
Let $v$ be a 1-periodic nonconstant real analytic function on $\mathbb{R}$. Then there exists some $\rho>0$ such that
\[
  v(x)=\sum_{k\in \mathbb{Z}}\hat v(k)e^{2\pi ikx},\ \ \mbox{with}\ \ |\hat v(k)|\lesssim e^{-\rho |k|}.
\]
Thus, there is a holomorphic extension $$v(z)=\sum_{k\in \mathbb{Z}}\hat v(k)e^{2\pi ikz}$$
to the strip $|\Im z|<\frac{\rho}{5}$, satisfying
\begin{equation}\label{supermum}
  |v(z)|\leq \sum_{k\in \mathbb{Z}}|\hat v(k)|e^{2\pi |k||\Im z|}<\sum_{k\in\mathbb{Z}}e^{-\rho |k|}e^{\rho |k|\frac{\pi}{5}}<C_v.
\end{equation}
And we need the following lemma from \cite{BG}, which only holds for the analytic functions on $\mathbb{T}$:
\begin{lemma}[Lemma 14.5 in \cite{BG}]
  For all $0<\delta <\rho$, there is  an $\epsilon$ such that
  \begin{equation}\label{es1}
    \inf_{E_1}\sup_{\frac{\delta}{2}<y<\delta}\inf_{x\in [0,1]}|v(x+iy)-E_1|>\epsilon.
  \end{equation}
\end{lemma}~\\

Now we begin to prove Theorem \ref{positive}. By Remark \ref{re} and (\ref{supermum}), we only need to assume that $|E|<(C_v+1)\lambda$ in the following paper.  Then $M_n(z,y,E)$ is analytic on $|\Im z|<\frac{\rho}{5}$ for fixed $y$ and $E$, with the norm $\|M_n(z,y,E)\|\leq (2C_v+2)^n\lambda^n$. Thus,  drop the fixed variables for convenience and define
 $$u_n(z):=\frac{1}{n}\log \|M_n(z)\|=\frac{1}{n}\log \|M_n(z,y,E)\|,$$
  which is a subharmonic function on $|\Im z|<\frac{\rho}{5}$, bounded by $\log [(2C_v+2)\lambda]<(1+\kappa) \log\lambda$ with  $\lambda>\lambda_1(v,\kappa)$.

Fix $0<\delta \ll \rho$ and $\epsilon$ satisfying Lemma 2.1. Define
\[\lambda_2=200\epsilon^{-\frac{2}{\kappa}}\]
and let $\lambda>\lambda_0(v,\kappa)=\max\{\lambda_1,\lambda_2\}$. Then, with fixed $E$, there is $\frac{\delta}{2}<y_0<\delta$ such that
\[
  \inf_{x\in[0,1]}\left |v(x+iy_0)-\frac{E}{\lambda}\right |>\epsilon,
\]
which implies that
\begin{equation}\label{es2}
  \inf_{x\in\mathbb{T}}|\lambda v(x+iy_0)-E|>\lambda \epsilon >200\epsilon^{-\frac{2}{\kappa}+1}>200.
\end{equation}
Let
\begin{equation}\label{setab}
  M_{n-1}(iy_0,E)\left (\begin{array}
  {cc} 1\\ 0
\end{array}\right )=\left (\begin{array}
  {cc} a_{n-1}\\ b_{n-1}
\end{array}\right ).
\end{equation}
Then
\begin{eqnarray}\label{anbn}\left (\begin{array}
  {cc} a_n\\ b_n
\end{array}\right )&=&\left (\begin{array}
  {cc} \lambda v\left (iy_0+\sum_{j=0}^{n-1}f(g^{j-1}(y))\right )-E &-1\\ 1&0\\
\end{array}\right )\left (\begin{array}
  {cc} a_{n-1}\\ b_{n-1}
\end{array}\right )\\
&=&\left (\begin{array}
  {cc} \left (\lambda v\left (iy_0+\sum_{j=0}^{n-1}f(g^{j-1}(y))\right )-E\right )a_{n-1}-b_{n-1}\\ a_{n-1}
\end{array}\right )\nonumber.\end{eqnarray}
Now we use the induction to show that
\[|a_n|\ge |b_n|,\ \ |a_n|>(\lambda\epsilon -1)|a_{n-1}|>(\lambda\epsilon-1)^n.\]
Due to (\ref{setab}) and (\ref{anbn}), it has that $a_0=1, b_0=0$ and 
\[|a_1|=|\lambda v(iy_0)-E|>\lambda\epsilon,\ |b_1|=1<\lambda\epsilon<|a_1|.\]
Let $|a_t|\ge |b_t|$ and $|a_t|>(\lambda\epsilon -1)|a_{t-1}|>(\lambda\epsilon-1)^t$. Then, we finish this induction by
\[|a_{t+1}|>\lambda \epsilon |a_{t}|-|b_{t}|>(\lambda\epsilon-1)|a_{t}|>(\lambda\epsilon-1)^t\ \mbox{and}\ |b_t|=|a_{t-1}|<|a_t|.\]
Thus, it implies that
\[
  \|M_n( iy_0, E)\|>\left|\left<M_n( iy_0, E)\left(\begin{array}
    {cc}1\\ 0\\
  \end{array}\right),\left(\begin{array}
    {cc}1\\ 0\\
  \end{array}\right)\right>\right|=|a_n|>(\lambda \epsilon-1)^n, 
\]and 
\[u_n(iy_0)>\log(\lambda \epsilon-1).\]

Write $\mathbb{H}=\{z:\Im z>0\}$ for the upper half-plane and $\mathbb{H}_s$ for the strip $\{z=x+iy:0<y<\frac{\rho}{5}\}$. Denote by $\mu(z,\mathcal{E},\mathbb{H})$ the harmonic measure of $\mathcal{E}$ at $z\in\mathbb{H}$ and $\mu_s(iy_0,\mathcal{E}_s,\mathbb{H}_s)$ the harmonic measure of $\mathcal{E}_s$ at $iy_0\in\mathbb{H}_s$, where $\mathcal{E}\subset\partial \mathbb{H}=\mathbb{R}$ and $\mathcal{E}_s\in \partial \mathbb{H}_s=\mathbb{R}\bigcup [y=\frac{\rho}{5}]$. Note that $\psi(z)=\exp\left (\frac{5\pi}{\rho}z\right )$ is a conformal map from $\mathbb{H}_s$ onto $\mathbb{H}$. Due to \cite{GM}, we have
\[
  \mu_s(iy_0,\mathcal{E}_s,\mathbb{H}_s)\equiv \mu(\psi(iy_0),\psi(\mathcal{E}_s),\mathbb{H}),
\]
and
\[
  \mu(z=x+iy,\mathcal{E},\mathbb{H})=\int_{\mathcal{E}}\frac{y}{(t-x)^2+y^2}\frac{dt}{\pi}.
\]
Thus
\[\mu_s[y=\frac{\rho}{5}]=\frac{5\pi y_0}{\pi\rho}<\frac{5\delta}{\rho}.\]
By the subharmonicity, it yields
\begin{eqnarray}
  \log (\lambda \epsilon-1)<u_n(iy_0)& \leq &\int_{[y=0]\bigcup [y=\frac{\rho}{5}]}u_n(z)\mu_s(dz)\nonumber\\
  &=&\int_{y=0}u_n(x)\mu_s(dx)+\int_{y=\frac{\rho}{5}}u_n(x+iy)\mu_s(dx)\nonumber\\
  &\leq &\int_{\mathbb{R}}u_n(x)\mu_s(dx)+\frac{5\delta}{\rho}\left [\sup_{y=\frac{\rho}{5}}u_n(x+iy)\right ]\nonumber\\
  &\leq &\int_{\mathbb{R}}u_n(x)\mu_s(dx)+\frac{\bar C\delta}{\rho}\log \lambda\nonumber.
\end{eqnarray}
So, by the setting of $\lambda_0$ and $\delta \ll \rho$, we have
\begin{eqnarray}
  \int_{\mathbb{R}}u_n(x)\mu_s(dx)&\ge& \log (\lambda\epsilon-1)-\frac{\bar C\delta}{\rho}\log \lambda\nonumber\\
  &\ge & \left (1-\frac{\bar C\delta}{\rho}\right )\log \lambda+\log \epsilon\nonumber\\
  &> &\left (1-\kappa\right )\log \lambda.\label{es-un}
\end{eqnarray}
Set
\[u^h_n(x)=u_n(x+h),\ \ h\in \mathbb{T}.\]
Then, due to (\ref{es2}), it is easy to see that (\ref{es-un}) also holds for $u^h_n(x)$. So, for any $h\in\mathbb{T}$, it has
\[\int_{\mathbb{R}}u_n(x+h)\mu_s(dx)> \left (1-\kappa\right )\log \lambda.\]
Integrating in $h\in \mathbb{T}$ implies that
\begin{eqnarray}
  L_n(y,E)=\int_0^1u_n(x+h)dh&\ge & \left (\int_{\mathbb{R}}\mu_s(dx)\right )\times\left ( \int_0^1u_n(x+h)dh\right )\nonumber \\
  &=&\int_0^1\int_{\mathbb{R}}u_n(x+h)\mu_s(dx)dh \nonumber\\
  &>&\left (1-\kappa\right )\log \lambda,\ \ \forall n\ge 0.\label{lny}
\end{eqnarray}
Thus
\[L_n(E)=\int_\mathbb{Y}L_n(y,E,\omega)dm(y)>\left (1-\kappa\right )\log \lambda,\ \ \forall n\ge 0,\]
which finishes this proof with $n\to +\infty$.

\section{Large Deviation Theorem}
From now on, we begin to consider the Schr\"odinger equations (\ref{schequ}).
For ease, we assume  $d=2$.  Then
\[v(\pi_1(T^n_\omega(x,y)))=v(x+ny+\frac{n(n-1)}{2}\omega),\]
\[ M_n(x,y,E,\omega)=\prod_{j=n}^1\left ( \begin{array}{cc}
 E-\lambda v(x+jy+\frac{j(j-1)}{2}\omega) & -1 \\
 1& 0 \\
  \end{array}\right ),\]and
  \[ L_n(E,\omega)=\frac{1}{n}\iint_{\mathbb{T}^2}\log \| M_n(x,y,  E,\omega)\|dxdy.\]
Recall that with fixed $y,E$ and $\omega$, $u_n(z)=\frac{1}{n}\log \|M_n(z,y,E,\omega)\|$ is a subharmonic function on $|\Im z|<\frac{\rho}{5}$ with the upper bound $\frac{101}{100}\log \lambda$($\kappa=\frac{1}{100}$). So we declare that the Fourier coefficient of $u_n(x)$ satisfies
\begin{equation}\label{fc}
  |\hat u_n(k)|\lesssim \frac{C}{|k|},\ \forall k\not=0.
\end{equation}Here we will use the following lemma(Lemma 2.2 in \cite{GS1}) to show that this constant $C$ depends only on $\lambda$ and $v$, but does not depend on $y$, $E$ or $\omega$:
\begin{lemma}
\label{lem:riesz} Let $u:\Omega\to \IR$ be a subharmonic function on
a domain $\Omega\subset\IC$. Suppose that $\partial \Omega$ consists
of finitely many piece-wise $C^1$ curves. There exists a positive
measure $\mu$ on~$\Omega$ such that for any $\Omega_1\Subset \Omega$
(i.e., $\Omega_1$ is a compactly contained subregion of~$\Omega$),
\[ u(z) = \int_{\Omega_1}
\log|z-\zeta|\,d\mu(\zeta) + h(z),
\]
where $h$ is harmonic on~$\Omega_1$ and $\mu$ is unique with this
property. Moreover, $\mu$ and $h$ satisfy the bounds \begin{eqnarray*}
  \mu(\Omega_1) &\le& C(\Omega,\Omega_1)\,(\sup_{\Omega} u - \sup_{\Omega_1} u), \\
\|h-\sup_{\Omega_1}u\|_{L^\infty(\Omega_2)} &\le&
C(\Omega,\Omega_1,\Omega_2)\,(\sup_{\Omega} u - \sup_{\Omega_1} u)
\end{eqnarray*} for any $\Omega_2\Subset\Omega_1$.
\end{lemma}
\noindent Thus,  there exists a constant $C=C(\lambda,v)$ such that for any $y$, $\omega$ and $E$,
\[\|\mu\|+\|h\|\leq C.\]
Then (\ref{fc}) holds by Corollary 4.7 in \cite{B}.~\\

Note that
\[u_n(\pi_1(T^j_\omega (x,y)))=u_n(x+jy+\frac{j(j-1)}{2}\omega)=L_n(y,E,\omega)+\sum_{k\not=0}\hat u(k)e^{ik(x+jy+\frac{j(j-1)}{2}\omega)}.\]
Then,
\begin{eqnarray*}\label{deq0}
  \left |\frac{1}{N}\sum_{j=1}^Nu_n(\pi_1(T^j_\omega (x,y)))-L_n(y,E,\omega)\right |
&=&\left |\frac{1}{N}\sum^N_{j=1}\sum_{k\in
\mathbb{Z}\backslash \{0\}}\hat u_n(k)e^{ik(x+jy+\frac{j(j-1)}{2}\omega)}\right |\\
&\leq &\frac{1}{N}\left | \sum_{0<|k|\leq K}\hat{u}_n(k)\sum^N_{j=1}e^{ik(jy+\frac{j(j-1)}{2}\omega)}\right |+\frac{1}{N}\left | \sum_{|k|> K}\hat{u}_n(k)\sum^N_{j=1}e^{ik(jy+\frac{j(j-1)}{2}\omega)}\right |\nonumber\\
&:=&(a)+(b)\nonumber
\end{eqnarray*}
Due to (\ref{fc}), it has
\begin{equation}\label{deqb}
  \|(b)\|^2_2\leq \sum_{|k|>K}\left |\hat u_n(k)\right |^2\leq C^2K^{-1}.
\end{equation}
To estimate (a), we will be using the following well-known method of Weyl-differencing, from Lemma 12 in \cite{K}:
\begin{lemma}[Weyl-differencing]\label{expsum}
Let $f(x)$ be a polynomial of degree $d\ge 2$:
\[f(x)=a_0+a_1x+\cdots+a_dx^d.\]
Then for any $k\ge 1$, it has
  \[
    \left |\sum_{x=1}^P e^{ i f(x)}\right |^{2^k}\leq 2^{2^k-1}P^{2^k-(k+1)}\sum_{y_1=0}^{P_1-1}\cdots\sum_{y_k=0}^{P_k-1}\left |\sum_{x=1}^{P_{k+1}}e^{ i\triangle_{y_1,\cdots,y_k}f(x)}\right |,
  \]
  where $P_1=P$ and under $\nu=1,2\cdots,k,$ quantities $P_{\nu+1}$ are determined by the equality $P_{\nu+1}=P_{\nu}-y_{\nu}$. Here $\triangle_{y_1} f(x)$ denotes  the  finite difference of a function $f(x)$ with an integer $y_1>0$:
\[\triangle_{y_1}f(x)=f(x+y_1)-f(x),\]
and when $k\ge 1$, the finite difference of the k-th order $\triangle_{y_1,\cdots,y_k} f(x)$ is determined with the help of the equality
\[\triangle_{y_1,\cdots,y_k} f(x)=\triangle_{y_k}\left [\triangle_{y_1,\cdots,y_{k-1}}f(x)\right ].\]
\end{lemma}
\noindent So, we have
\[
  \left |\sum^N_{j=1}e^{ik(jy+\frac{j(j-1)}{2}\omega)}\right |^2\leq 2N+2\sum_{m=1}^{N-1}\min \left(N-m,2\|km\omega\|^{-1}\right )\leq C_1 \sum_{m=1}^{N-1}\min \left(N,2\|km\omega\|^{-1}\right ).\nonumber
\]
And due to the Cauchy inequality, it has
\begin{eqnarray}
  |(a)|^2&\leq& N^{-2}\sum_{0<|k|\leq K}\left |\hat{u}(k)\right |^2\sum_{0<|k|\leq K}\left (\sum^N_{j=1}e^{ik(jy+\frac{j(j-1)}{2}\omega)}\right )^2\leq C_2N^{-2}\sum_{k=1}^K\sum_{m=1}^{N-1}\min \left(N,2\|km\omega\|^{-1}\right )\nonumber\\
  &\leq &C_{\epsilon}N^{-2}(KN)^{\epsilon}\sum_{k=1}^{KN}\min \left(N,2\|k\omega\|^{-1}\right ),\label{epsilonap}
\end{eqnarray}
where the arbitrary small positive parameter $\epsilon0$ and (\ref{epsilonap}) come from the following lemma, which is Lemma 13 in \cite{K}:
\begin{lemma}\label{pronum}
  Let $M$ and $m_1,m_2,\cdots,m_n$ be positive integers. Denote by $\tau_n(M)$ the number of solutions of the equation $m_1\cdots m_n=\lambda$. Then under any $\epsilon (0<\epsilon\leq 1)$ we have
  \[ \tau_n(M)\leq C_n(\epsilon)M^{\epsilon},\]
  where the constant $C_n(\epsilon)$ depends on $n$ and $\epsilon$ only.
\end{lemma}
\noindent By Dirichlet's principle there is an integer $1\leq q \leq N$ and an integer $p$ so that $gcd(p,q)=1$ and $\left |\omega-\frac{p}{q}\right |\leq \frac{1}{qN}$. Thus duo to the definition of Diophantine number, one has
\begin{equation}\label{dp}
  N\ge q\ge c_{\omega}\frac{N}{(\log N)^2}.
\end{equation}
Combined with (\ref{dp}), the following lemma, which is also from \cite{K}(Lemma 14), will help us evaluate (\ref{epsilonap}):
\begin{lemma}\label{de}
  Let $P\ge 2$ and
  \[\omega=\frac{p}{q}+\frac{\theta}{q^2},\ \ (p,q)=1,\ \ |\theta|\leq 1.\]
  Then under any positive integer $Q$ and an arbitrary real $\beta$ we have
  \[
    \sum_{x=1}^Q\min \left (P,\frac{1}{\|\omega x+\beta\|}\right )\leq 4\left (1+\frac{Q}{q}\right )(P+q\log P).
  \]
\end{lemma}
\noindent Thus,
\begin{eqnarray}\label{sum1}
  |(a)|^2&\leq &C_{\epsilon}N^{-2}(KN)^{\epsilon}\sum_{k=1}^{KN}\min \left(N,2\|k\omega\|^{-1}\right )\\
  &\leq & 4C_{\epsilon}N^{-2}(KN)^{\epsilon}\left (1+\frac{KN}{q}\right )(N+q\log N)\nonumber\\
  &\leq&4C_{\epsilon}N^{-2}(KN)^{\epsilon}\left (N+N\log N+\frac{1}{c_\omega}KN(\log N)^2+KN\log N\right )\nonumber\\
  &\leq & N^{2\epsilon-1}K^{1+2\epsilon},\nonumber
\end{eqnarray}for $N>N_0(\epsilon,\omega)$.
Note that
\[|u_n(x)-u_n(\pi_1(T_{\omega}(x,y)))|<\frac{C}{n}.\]
Then,
\[
  \left |u_n(x)-\frac{1}{N}\sum_{j=1}^Nu_n(\pi_1(T_{\omega}^j(x,y)))\right |<\frac{CN}{n}.
\]
Above all, we have
\begin{equation}\label{es5}
 \left |u_n(x)-L_n(y,E,\omega)\right |\leq (a)+(b)+\frac{CN}{n},
 \end{equation}
 with the estimations  (\ref{deqb}) and (\ref{sum1}).~\\

 To improve (\ref{es5}), we will apply the following lemma proved in \cite{BGS}, which gives the evaluation of the BMO norm for subharmonic functions:
\begin{lemma}[Lemma 2.3 in \cite{BGS}]\label{BMO}
Suppose u is subharmonic on $\mathscr{A}_{\rho}$, with $\sup_{\mathscr{A}_{\rho}}|u|\leq n$. Furthermore, assume that $u=u_0+u_1$, where
\begin{equation}\label{setting1}
\|u_0-<u_0>\|_{L^{\infty}(\mathbb{T})}\leq \epsilon_0\ \ \mbox{and}\ \ \|u_1\|_{L^1(\mathbb{T})}\leq \epsilon_1.
\end{equation}
Then for some constant $C_{\rho}$ depending only on $\rho$,
\[\|u\|_{BMO(\mathbb{T})}\leq C_{\rho}\left (\epsilon_0\log \left(\frac{n}{\epsilon_1}\right )+\sqrt{n\epsilon_1}\right).\]
\end{lemma}
\noindent Now, if we choose $u(z)=\frac{n}{C}u_n(z)$, then $ u(z)$ is subharmonic on $\mathscr{A}_{\rho}$, with $\sup_{\mathscr{A}_{\rho}}| u|\leq n$. Let $K=n^{\frac{1}{4}}$, $N=n^{\frac{9}{10}}$ and $\epsilon=\frac{1}{5}$. By (\ref{es5}), it has
\begin{equation}\label{bmo1}
  \left | u(x)-< u(\cdot)>\right |\leq 3n^{\frac{9}{10}},
\end{equation}
up to a set of measure less than $n^{-\frac{1}{40}}$. Define $\mathscr{B}$ to be this exceptional set and let
\[ u(x)-< u(\cdot)>=u_0+u_1,\]
where $u_0=0$ on $\mathscr{B}$ and $u_1=0$ on $\mathbb{T}\backslash \mathscr{B}$. Then,
\[ \|u_0-<u_0>\|_{L^{\infty}(\mathbb{T})}\leq 3n^{\frac{9}{10}}\ \ \mbox{and}\ \ \|u_1\|_{L^1(\mathbb{T})}\leq n^{\frac{39}{40}} .\]
Due to Lemma \ref{BMO}, it yields
\[\| u\|_{BMO(\mathbb{T})}\leq C_{\rho}n^{\frac{79}{80}} .\]
 Recall the following John-Nirenberg inequality (\cite{S}):
\[meas\{x\in\mathbb{T}: | u(x)-< u>|>\gamma \}\leq C\exp \left (-\frac{c\gamma}{\| u\|_{BMO}}\right ),\]
with the absolute constants $C$ and $c$. Let $\gamma= n\frac{1}{50}\log \lambda$. Then there exists $n_0(\omega,\lambda,v)$ such that for any  $y\in\mathbb{T}$, $|E|<(C_v+1)\lambda$ and $n>n_0$,
   \begin{equation}\label{eqldt}
     meas\{x\in\mathbb{T}: |\frac{1}{n}\log\|M_n(x,y,E,\omega)\|-L_n(y,E,\omega)|>\frac{1}{50}\log \lambda \}\leq \exp C (-cn^{\frac{1}{80}}\log \lambda ).
   \end{equation}
We emphasize again that $C$ and $c$ are the absolute constants, not depend on $y$. Thus, combining Theorem \ref{positive}, (\ref{lny}) and (\ref{eqldt}), we have the following Large Deviation Theorem:
\begin{lemma}[Large Deviation Theorem]\label{ldt}
  Let $\lambda_0(v)$ be as in Theorem \ref{positive}  with $\kappa=\frac{1}{100}$. Assume $\lambda>\lambda_0$ and $n>n_0(\omega,\lambda,v)$. Then
   \[meas\{(x,y)\in\mathbb{T}^2: |\frac{1}{n}\log\|M_n(x,y,E,\omega)\|-L_n(E,\omega)|>\frac{1}{40}L_n(E,\omega)\}\leq C\exp  (-cn^{\frac{1}{80}}\log \lambda ) .\]
   \end{lemma}

\section{Proof of Theorem \ref{skt}}
Actually, the non-perturbation of Theorem \ref{skt} comes directly from the non-perturbation of Large Deviation Theorem, Lemma \ref{ldt}. Of course, it also needs some other methods, such as avalanche principle, Green function estimate, semi-algebraic set theory, parameterizations and so on. But these methods had been developed and worked well for the shew-shift Schr\"odinger equations in \cite{BGS} and \cite{K2}, and we can apply them directly to the same equations in our paper. Thus in this section, for readers' ease, we only give the main idea of the proofs and point out how our LDT works in them. The readers can see the details in the above two references.

\subsection{Week-H\"older Continuity of the Lyapunov exponent}
  First, we will use the induction to show that
  \begin{equation}\label{inf-est}
    \left |L(E)-2L_{2n_0}(E)+L_{n_0}(E)\right |\leq \exp(-cn_0^{\frac{1}{1000}}\log \lambda).
  \end{equation}
  Indeed, by Lemma \ref{ldt} and the avalanche principle, we can construct a series of $\{n_j\}_{j=0}^{\infty}$, which satisfies $$n_{j+1}\simeq\exp(cn_j^{\frac{1}{200}}\log \lambda),\ \ j=0,1,\cdots,$$ such that
  \[ \left |L_{n_{j+1}}(E)-2L_{2n_j}(E)+L_{n_j}(E)\right |\leq \exp(-cn_j^{\frac{1}{200}}\log \lambda),\]
  and
  \[ \left |L_{2n_{j+1}}(E)-2L_{2n_j}(E)+L_{n_j}(E)\right |\leq \exp(-cn_j^{\frac{1}{200}}\log \lambda).\]Thus, we have
  \[ \left |L_{2n_{j+1}}(E)-L_{n_{j+1}}(E)\right |\leq 2 \exp(-cn_j^{\frac{1}{200}}\log \lambda),\]
  and
  \[ \left |L_{n_{j+1}}(E)-L_{n_{j}}(E)\right |\leq 5 \exp(-cn_{j-1}^{\frac{1}{200}}\log \lambda).\]Then
  \begin{eqnarray}
    (\ref{inf-est})&\leq & \sum_{j\geq 1}|L_{n_{j+1}}(E)-L_{n_{j}}(E)|+|L_{n_1}(E)+L_{n_0}(E)-2L_{2n_0}(E)|\nonumber\\
                   &\leq & \sum_{j\geq 1}5 \exp(-cn_{j-1}^{\frac{1}{200}}\log \lambda)+ \exp(-cn_{0}^{\frac{1}{200}}\log \lambda)\nonumber\\
                   &\leq &  \exp(-cn_0^{\frac{1}{1000}}\log \lambda).\nonumber
  \end{eqnarray}

  Second, using Trotter's formula, we have
\[\left |\log\|M_{n}(\un x,E')\|-\log \|M_{n}(\un x,E)\| \right |\leq\big |\|M_{n}(\un x,E')\|-\|M_{n}(\un x,E)\|\big |
\leq n(C_v \lambda)^{n-1}|E'-E|.
\]It implies that
\[\left |L_n(E)-L_n(E')\right |<(C_v \lambda)^{n-1}|E'-E|.\]

Finally,
\begin{eqnarray}
  \left |L(E)-L(E')\right |&\leq &
\left |L(E)-2L_{2n_0}(E)+L_{n_0}(E)\right |+\left |L(E')-2L_{2n_0}(E')+L_{n_0}(E')\right |\nonumber\\
&\ &\ \ +|L_{n_9}(E)-L_{n_0}(E')|+2|L_{2n_0}(E)-L_{2n_0}(E')|\nonumber\\
&< &2\exp(-cn_0^{\frac{1}{1000}}\log \lambda)+3(C_v \lambda)^{n_0-1}|E'-E|\leq \exp\left (-c\left |\log |E-E'|\right |^{\tau}\right ),\nonumber
\end{eqnarray}
if one sets $|E-E'|=\exp(-n_0\log \lambda)$ with large $n_0$.

\subsection{Non-perturbative Anderson Localization}
Define
\begin{eqnarray*}
  S^N_{\un x,\omega,\lambda v}&=&R^{[1,N]}S_{\un x,\omega,\lambda v}R^{[1,N]}\\
  &=&\left(\begin{array}
  {cccccc}
  \lambda v\left(\pi_1\left(T_{\omega}(x,y)\right)\right)&-1&0&\cdots&\cdots&0\\
  -1&\lambda v\left(\pi_1\left(T^2_{\omega}(x,y)\right)\right)&-1&0&\cdots&0\\
  \vdots&\vdots&\vdots&\vdots&\vdots&\vdots\\
  0&\cdots&0&-1&\lambda v\left(\pi_1\left(T^{N-1}_{\omega}(x,y)\right)\right)&-1\\
  0&\cdots&\cdots&0&-1&\lambda v\left(\pi_1\left(T^N_{\omega}(x,y)\right)\right)\\
\end{array}\right),
\end{eqnarray*}
where $R^{[1,N]}$ is the coordinate restriction matrix to $[1,N]\subset \mathbb{Z}$. Then by Cramer's rule, one has for $1\leq N_1\leq N_2\leq N$,
\begin{eqnarray}\label{upest1}
  \left |G_N(\un x,E,\omega)(N_1,N_2)\right |&=& \left |(S^N_{\un x,\omega,\lambda v}-E)^{-1}(N_1,N_2)\right |\nonumber\\
  &=&\frac{\left |\det [S^{N_1-1}_{\un x,\omega,\lambda v}-E)]\right |\left |\det[S^{N-N_2}_{T^{N_2}\un x,\omega,\lambda v}-E)]\right | }{\left |\det[S^{N}_{\un x,\omega,\lambda v}-E)]\right |}\nonumber\\
  &\leq &\frac{\left \|M_{N_1}(\un x,E,\omega)\right \|\left \|M_{N-N_2}(T^{N_2}\un x,E,\omega)\right \| }{\left |\det[S^{N}_{\un x,\omega,\lambda v}-E)]\right |},
\end{eqnarray}
where (\ref{upest1}) comes from the following relationship between $M_n$ and determinants:
\begin{equation}\label{re-det}
 M_n(\un x,E,\omega)=\left [\begin{array}
  {cc} \det [S^{n}_{\un x,\omega,\lambda v}-E)]& -\det [S^{n-1}_{T\un x,\omega,\lambda v}-E)]\\
  \det [S^{n-1}_{\un x,\omega,\lambda v}-E)]&-\det [S^{n-2}_{T\un x,\omega,\lambda v}-E)]
\end{array}
\right ].
\end{equation}
Note that if $\lambda>\lambda_0$, then $\|M_n(\un x,E,\omega)\|\leq \exp(\frac{101}{100}n\log \lambda)$ for any $|E|\leq C_v\lambda+1$. Thus
\begin{equation}\label{upest2}
  (\ref{upest1})\leq \frac{\exp\left (\frac{101}{100} [N-|N_1-N_2|]\log \lambda\right ) }{\left |\det[S^{N}_{\un x,\omega,\lambda v}-E)]\right |}.
\end{equation}
Returning to (\ref{re-det}), if we allow replacement of $N$ by $N-1$ or $N-2$ and $\un x$ by $T \un x$, may replace the denominator in (\ref{upest2}). Then, by Theorem \ref{positive} and Lemma \ref{ldt}, we have
\[
\left |G_\Lambda(\un x,E,\omega)(N_1,N_2)\right |\leq \exp\left (-\frac{101}{100} |N_1-N_2|\log \lambda+\frac{1}{20} N\log\lambda\right ),
\]up to a set $\Omega_N(E)$ satisfying $meas [\Omega_N(E)]<\exp  (-cN^{\frac{1}{80}}\log \lambda )$, where $\Lambda$ is one of the intervals
\[[1,N],\ \ [1,N-1],\ \ [2,N],\ \ [2,N-1].\]
So, if $|N_1-N_2|\ge \frac{1}{10}N$, then
\begin{equation}\label{goodgreen}
  \left |G_\Lambda(\un x,E,\omega)(N_1,N_2)\right |\leq \exp\left (-\frac{1}{20}  N\log \lambda\right ).
\end{equation}
It is easy to see that (\ref{goodgreen}) also holds  with $|N_1-N_2|\ge \frac{1}{10}N$ if we redefine $\Lambda$ to be one of the following intervals
\[[-N,N],\ \ [-N,N-1],\ \ [-N+1,N],\ \ [-N+1,N-1].\]
Define $\Omega=\Omega(E)$ to be the exception set for (\ref{goodgreen}) with measure $<\exp  (-cN^{\frac{1}{80}}\log \lambda)$. Here we need to emphasize again that this set depends on the energy $E$.~\\

Fix $\un x_0\in \mathbb{T}^2$ and consider the orbit $\{T_\omega^j(\un x_0): |j|\le N_1\}$ where $N_1=N^{C'}$. Then by Lemma 15.21 in \cite{B}, we have
\[
  \sharp \{|j|\le N_1:T_\omega^j(\un x_0)\in \Omega\}<N_1^{\frac{19}{20}}.
\]
Thus, except for at most $N_1^{\frac{19}{20}}$ values of $|j|<N_1$,
\[\left |G_\Lambda(T_\omega^j(\un x_0),E,\omega)(N_1,N_2)\right |\leq \exp\left (-\frac{1}{20}  N\log \lambda\right ),\]
if $|N_1-N_2|\ge \frac{1}{10}N$. Note that  the value of any formal solution $\phi$ of the equation $S_{\un x_0,\omega}\phi=E\phi$ at  a point $n\in [a,b]\subset \mathbb{Z}$ can be reconstructed from the boundary values via
\begin{equation}\label{greenexpre}
  \phi(n)=G_{[a,b]}(\un x_0,E,\omega)(n,a)\phi(a-1)+G_{[a,b]}(\un x_0,E,\omega)(n,b)\phi(b+1).
\end{equation}
Thus, if $E$ is a generalized eigenvalue and $\phi$ is the corresponding generalized eigenfunction, then by (\ref{greenexpre}),
\begin{equation}\label{goodj}
 |\phi(j)|<\exp\left (-\frac{1}{100}N\log \lambda \right )
\end{equation}
holds for all $|j|<N_1$ except $N_1^{\frac{19}{20}}$ many. Let both of $j_0$ and  $-j_0$  satisfy (\ref{goodj}). By (\ref{greenexpre}) again, we have
\[1\leq \|G_{[-j_0+1,j_0-1]}(\un x_0,E,\omega)\|\bigg (|\phi(j_0)|+|\phi(-j_0)|\bigg ).\]
Thus
\[\|G_{[-j_0+1,j_0-1]}(\un x_0,E,\omega)\|\ge \exp\left (\frac{1}{100}N\log \lambda \right ),\]
which also means
\begin{equation}\label{closed-e}
\mbox{dist} \left (E, \mbox{Spec}\  S_{\un x_0,\omega,\lambda v}^{[-j_0+1,j_0-1]}\right )<\exp\left (-\frac{1}{100}N\log \lambda \right ).
\end{equation}
Define $\mathcal{E}_{\omega}=\bigcup_{|j|\le N_1} \mbox{Spec}\  S_{\un x_0,\omega,\lambda v}^{[-j_0+1,j_0-1]}$. So, by (\ref{closed-e}), if $\un x\not \in \bigcup_{E'\in \mathcal{E}_{\omega}}\Omega(E')$ and $E$ is a generalized eigenvalue of $S_{\un x_0,\omega,\lambda v}$, then we have 
\[\left |G_\Lambda(\un x,E,\omega)(N_1,N_2)\right |\leq \exp\left (-\frac{1}{150}  N\log \lambda\right )\]
 with $|N_1-N_2|\ge \frac{1}{10}N$.~\\

At last, we  show that the measure of the set of the $\omega$, which makes $T_{\omega}^j (\un x)$ belongs to this "bad set" $\bigcup_{E'\in \mathcal{E}_{\omega}}\Omega(E')$ for some $|j|\sim M= N^{C"}(C"\gg C'> 1)$, is zero. Fortunately, the semi-algebraic set theory can help us prove the following lemma:
\begin{lemma}[Lemma 15.26 in \cite{B}]
Let $S\in\mathbb{T}^3$ be a semi-algebraic set of degree $B$, s.t.
\[meas[S]<\exp(-B^{\sigma}),\ \mbox{for } \sigma>0.\]
Let $M$ be an integer satisfying
\[\log\log M\ll \log B \ll \log M.\]
Thus, for any fixed $\un x\in\mathbb{T}^2$
\[
  meas [\omega\in \mathbb{T}|\left (\omega, T_{\omega}^j(\un x)\right )\in S \ \mbox{for some }j\sim M]<M^{-c}
\]
for some $c>0$.
\end{lemma}
\noindent Let $\mathcal{R}_M$ denote the above $\omega$-set. Then, $ |\phi(n)|<\exp\left (-\frac{1}{200}n\log \lambda \right ) $  for any $M^{\frac{1}{2}}<|n|<2M$ and $\omega\not\in \mathcal{R}_M$. Set
\[\mathcal{R}=\bigcup_{N}\bigcap_{M>N}\mathcal{R}_M.\]
We have $meas[\mathcal{R}]=0$ and for $\omega \in DN\backslash \mathcal{R}$, the non-perturbative Anderson localization holds.~\\~\\

\subsection{Intervals In The Spectrum}By minimality of the skew-shift, $\sigma(S_{\un x,\omega,\lambda v})=\sigma(S_{\un y,\omega,\lambda v})$ for any $\un x,\ \un y\in \mathbb{T}^2$. Thus, the following theorem implies our Statement (I) directly:
\begin{theorem}\label{si}
  Let $\lambda>\max\{\lambda_0,\lambda_1\}$. Then for any $E\in\lambda \mathcal{E}_\delta$, there exists $\un x\in\mathbb{T}^2$ such that $E$ is an eigenvalue of $S_{\un x,\omega,\lambda v}$.
\end{theorem}

The proof of this theorem is concerned with the following eisolated eigenvalues and parametrization:
\begin{defi}
  Let  A be a self-adjoint operator, $\epsilon>0$ and $E\in\mathbb{R}$. E is an $\epsilon-$isolated eigenvalue of A, if
  \[\sigma(A)\bigcap [E-\epsilon,E+\epsilon]=\{E\}\]
  and $E$ is simple.
\end{defi}
\begin{defi}
  Let $\zeta:\mathcal{X}\subset \mathbb{T}\to\mathbb{T}$ be a continuously differential function, $\epsilon>0$, $L\in (0,\frac{1}{3})$ and $M\ge 0$. A pair $(\zeta, \mathcal{X})$ is called to be an $(\epsilon,L)$-parametrization of the eigenvalue $E_0$ of $S_{\cdot, \omega,\lambda v}^{[-M,M]}$, if
  \begin{enumerate}
    \item[{\rm{(1)}}]For any $x\in \cx$, we have that $E_0$ is an $\epsilon$-isolated eigenvalue of $S_{(x,\zeta(x)), \omega,\lambda v}^{[-M,M]}$;
    \item[{\rm{(2)}}]$\me[\cx]\ge\frac{1}{\sqrt{\max\{M,1\}}}$;
    \item[{\rm{(3)}}]$\|\zeta'\|_{L^\infty(\mathbb{T})}\le L$.
  \end{enumerate}
\end{defi}
\noindent
Now what we want is to use the induction to structure a sequence of $\{(\zeta_j, \mathcal{X}_j)\}_{j=1}^\infty$ which is $(\epsilon_j,L_j)$-parametrization of the eigenvalue $E_0$ of $S_{\cdot, \omega,\lambda v}^{[-M_j,M_j]}$, satisfying
\begin{enumerate}
    \item[{\rm{(i)}}]there exists some $x_{\infty}\in \bigcap_{j=1}^\infty \cx_j\not =\emptyset$;
    \item[{\rm{(ii)}}]$y_j=\zeta_j(x_{\infty})\to y_\infty$ in $l^2(\mathbb{Z})$;
    \item[{\rm{(iii)}}]the eigenfunctions $\psi_j$ of $S_{(x_{\infty},y_j), \omega,\lambda v}^{[-M_j,M_j]}$ corresponding to $E_0$ form a Cauchy sequence;
    \item[{\rm{(iv)}}]$M_j\to \infty$.
  \end{enumerate}
If we have (i)-(iv),  then by the continuity, it has
\[S_{(x_{\infty},y_\infty), \omega,\lambda v}\psi_\infty=E_0\psi_\infty,\]
where $\psi_\infty=\lim_{j\to\infty}\psi_j$. Then, we get Theorem \ref{si}.~\\

To obtain this sequence by the induction, we need to give the suitable $(\zeta_1, \mathcal{X}_1)$ and show how to structure $(\zeta_{j+1}, \mathcal{X}_{j+1})$ from $(\zeta_j, \mathcal{X}_j)$. Before this construction, we first make the following definition about the desired  communications   between these parameterizations:
\begin{defi}
  Let $(\zeta_{j}, \mathcal{X}_{j})$ be $(\epsilon_j,L_j)$-parametrization of the eigenvalue $E_0$ of $S_{\cdot, \omega,\lambda v}^{[-M_j,M_j]}$ for $j=1,2$.  We say $(\zeta_{2}, \mathcal{X}_{2})$ is a $\delta$-extension of $(\zeta_{1}, \mathcal{X}_{1})$, if
  \begin{enumerate}
    \item[{\rm{(1)}}]$\cx_2\subset \cx_1$;
    \item[{\rm{(2)}}]$\epsilon_2<\epsilon_1$, $M_2> M_1$;
    \item[{\rm{(3)}}]$L_2\leq L_1+\delta$ and $\|\zeta_1-\zeta_2\|_{L^\infty(\cx_2)}\le \delta$;
    \item[{\rm{(4)}}]Let $x\in\cx_2$ and $\psi_j\in l^2([-M_j,M_j])$ normalized eigenfunctions of $S_{\cdot, \omega,\lambda v}^{[-M_j,M_j]}$ corresponding to the eigenvalue $E_0$. We have for some $|a|=1$
        \[\left (\sum_{n\in\mathbb{Z}}(1+n^2)|\psi_1(n)-a\psi_2(n)|^2\right )^{\frac{1}{2}}\leq \delta.\]
  \end{enumerate}
\end{defi}
\noindent It is easy to check that if for any $j\ge 1$, $(\zeta_{j+1}, \mathcal{X}_{j+1})$ is a $\delta_{j+1}$-extension of $(\zeta_{j}, \mathcal{X}_{j})$ with $\delta_{j+1}\ll \delta_{j}\le \frac{1}{3}$, then the sequence $\{(\zeta_j, \mathcal{X}_j)\}_{j=1}^\infty$ satisfies (i)-(iv).
Now we can start the construction. The suitable $(\zeta_1, \mathcal{X}_1)$ comes directly from the following lemma:
\begin{lemma}[Theorem 3.6 in \cite{K2}]\label{fs}
  Let $M\ge 1$. Then there exists $\bar \lambda_0=\bar \lambda_0(M,v,\delta)>0$ such that for $\lambda>\bar \lambda_0$ and $E_0\in\lambda \mathcal{E}_\delta$, there exists a $\lambda^{-\frac{1}{500}}$--parametrization $(\zeta_1, \mathcal{X}_1)$ at scale $M$ that $\lambda^{-\frac{1}{10}}$--extension of $(\zeta_{0}, \mathcal{X}_{0})$.
\end{lemma}
\noindent Here $(\zeta_{0}, \mathcal{X}_{0})$ is the simplest one but not good enough:
\[M_0=0,\ \ \  S_{(x,y), \omega,\lambda v}^{[-0,0]}=\lambda v(y),\ \ \  \cx_0=\mathbb{T},\ \ \   \zeta_0(x)=y_0,\ \forall x\in \mathbb{T},\ \ \  E_0=\lambda v(y_0).\]
Moreover, the following theorem  tells us how to structure $(\zeta_{j+1}, \mathcal{X}_{j+1})$ from $(\zeta_{j}, \mathcal{X}_{j})$:
\begin{lemma}\label{ind}
  Let $M$ be large enough and $\lambda>\lambda_0$. Furthermore, assume for $\epsilon=\exp(-M^{\frac{1}{50}})$ that $(\zeta, \mathcal{X})$ is an $(\epsilon,L)$-parametrization of  the eigenvalue $E_0=\lambda v(y_0)$ of $S_{\cdot, \omega,\lambda v}^{[-M,M]}$ that $\frac{C_2^5}{2C_1}$-extends $(\zeta_{0}, \mathcal{X}_{0})$, where $C_1=10\|v'\|_{L^\infty(\mathbb{T})}$, $C_2=\frac{1}{10}\max\{|v'(y_0)|,1\}$ and $L+\epsilon\le \frac{1}{3}$. Define $R=\exp(M^{\frac{1}{1000}})$. Then there exists $( \zeta',  \mathcal{X}')$ such that $( \zeta', \mathcal{X}')$ is a $\left (\frac{1}{1000}\epsilon,L'\right )$-parametrization of  the eigenvalue $E_0$ of $S_{\cdot, \omega,\lambda v}^{[-R,R]}$ with $ L'=L+\epsilon$ and for $\eta=\exp(-\frac{1}{100}M)$, $( \zeta',  \mathcal{X}')$ is an $\eta$-extension of $(\zeta, \mathcal{X})$ to scale $R$.
\end{lemma}
\noindent  The proofs of the above two lemmas can be found in \cite{K2}. The first one is just Theorem 3.6 in that literature, whose proof is in Section 4, 9 and 10. The second one is a modified version of Kr\"uger's Theorem 3.5, as our non-perturbative Large Deviation Theorem replaces the perturbative one he used. Then, the other part of the proof can copy Section 5-8 in \cite{K2} directly.~\\

To end this paper, we give the induction. Let $M_1=M$. Then Lemma \ref{fs} tells us that there exists a $\lambda^{-\frac{1}{500}}$--parametrization $(\zeta_1, \mathcal{X}_1)$ at scale $M$ that $\lambda^{-\frac{1}{10}}$--extension of $(\zeta_{0}, \mathcal{X}_{0})$. If $\lambda$ is big enough to satisfy
\[\lambda^{-\frac{1}{500}}\leq \epsilon_1\ \mbox{and}\  \lambda^{-\frac{3}{2}}\leq \frac{C_2^5}{10C_1},\]
then the assumptions of Lemma \ref{ind} hold for $(\zeta_1, \mathcal{X}_1)$. Define the sequences
\[M_{j+1}\sim \exp(M_j^{\frac{1}{1000}}),\ \ \ \epsilon_j=\exp(-M_j^{\frac{1}{50}}),\ \ \ L_{j+1}=L_j+\epsilon_j,\ \ j=1,2,\cdots.\]
From Lemma \ref{ind}, we have that $(\zeta_{2}, \mathcal{X}_{2})$ is an $(\epsilon_2,L_2)$-parametrization of  the eigenvalue $E_0$ of $S_{\cdot, \omega,\lambda v}^{[-M_2,M_2]}$. Note that $(\zeta_{2}, \mathcal{X}_{2})$ is an $\exp(-\frac{1}{100}M_1)$-extension of $(\zeta_1, \mathcal{X}_1)$, also a $\frac{C_2^5}{2C_1}$-extension of  $(\zeta_{0}, \mathcal{X}_{0})$. Thus, $(\zeta_{2}, \mathcal{X}_{2})$ satisfies all the assumptions of Lemma \ref{ind}, which makes the induction work. In summary, Statement (I) is proved.

\end{document}